\documentclass[10pt]{amsart}
\newcommand{\ra}{\rightarrow}
\newcommand{\R}{\mathbb{R}}
\newcommand{\T}{\mathbb{T}}
\newcommand{\Z}{\mathbb{Z}}

\newcommand{\ep}{\epsilon}

\newcommand{\8}{\infty}

\newcommand{\ee}{\end{eqnarray}}
\newcommand{\be}{\begin{eqnarray}}
\newcommand{\dfr}{\mbox{{\rm Diff}$_\mu^r(M)$}}

\newtheorem{thm}{Theorem}[section]
\newtheorem{lem}[thm]{Lemma}

\newtheorem{conj}{Conjecture}
\begin{document}

\title{Area-Preserving Surface Diffeomorphisms}
\author{Zhihong Xia}
\address{Department of Mathematics \\ Northwestern University \\
Evanston, Illinois 60208}
\thanks{Research supported in part by National Science Foundation.}
\date{Februry 10, 2005. revised version}
\email{xia@math.northwestern.edu}

\begin{abstract}
  We prove some generic properties for $C^r$, $r=1, 2, \ldots,
  \infty$, area-preserving diffeomorphism on compact surfaces. The
  main result is that the union of the stable (or unstable) manifolds
  of hyperbolic periodic points are dense in the surface. This extends
  the result of Franks and Le Calvez \cite{FL03} on $S^2$ to general
  surfaces. The proof uses the theory of prime ends and Lefschetz
  fixed point theorem.
\end{abstract}

\maketitle

\section{Introduction and statement of main results}

Let $\dfr$ be the set of all $C^r$, $r=1, 2, \ldots, \infty$
diffeomorphisms of compact orientable surface $M$ that preserves a
smooth area element $\mu$.  A property for $C^r$ area-preserving
diffeomorphisms of $M$ is said to be {\em generic}\/, if there is a
residual subset $R \subset \dfr$ such that the property holds for all
$f \in R$. The following are some examples of known generic
properties:

\begin{itemize}
\item[G1] All the periodic points are either elliptic or
  hyperbolic. This is proved by Robinson \cite{Robinson70}.

\item[G2] All elliptic periodic points are Moser stable. Moser stable
  means that the normal form at the elliptic periodic points are
  non-degenerate in the sense of KAM theory (cf. Siegel \& Moser
  \cite{SM71}). This implies that there are invariant curves surrounding
  each elliptic periodic point. This condition requires that the map
  is sufficiently smooth, say $r \geq 16$.

\item[G3] For any two hyperbolic periodic points $p$, $q$, the
intersections of the stable manifold $W^s(p)$ and the unstable
manifold $W^u(q)$ are transversal. This is again proved by Robinson
\cite{Robinson70}. G1 and G3 together are often referred as the
Kupka-Smale condition for area-preserving diffeomorphisms.

\item[G4] For any hyperbolic periodic point $p$, let
$\Gamma_1$ and $\Gamma_2$ be any two branches of the stable manifold
  or unstable manifold of $p$, then $\overline{\Gamma_1}=
  \overline{\Gamma_2}$. For $S^2$, this is proved by Mather
  \cite{Mather82a}. The general case is proved by Oliveira
  \cite{Oliveira87}.

\item[G5] If $M$ is a two-sphere $S^2$ or a two torus $T^2$, let
$\Gamma_1$ be a branch of the stable manifold of a hyperbolic periodic
point of $p$ and $\Gamma_2$ be a branch of the unstable manifold of
$p$, then $\Gamma_1 \cap \Gamma_2 \neq \emptyset$. For $S^2$, this is
proved by Robinson \cite{Robinson73} and Pixton \cite{Pixton82}. For
$\T^2$, it is proved by Oliveira \cite{Oliveira87}. For general
surfaces, it is proved by Oliveira \cite{Oliveira00} for most homotopy
classes of maps. For $C^1$ case and any compact manifold of arbitrary
dimension, the result is proved by Takens \cite{Takens72} (see also
Xia \cite{Xi96a} for a stronger result).

\end{itemize}

Franks \& Le Calvez \cite{FL03} recently showed another remarkable
$C^r$ generic property for area-preserving diffeomorphisms on $S^2$. 
Their results state that the stable and unstable manifolds of
hyperbolic periodic points are dense in $S^2$ generically. In this
paper, we extend this result to arbitrary orientable surfaces and to
arbitrary homotopy classes. More precisely, we state our main result.

\begin{thm}
Let $M$ be a compact orientable surface. There exists a residual
subset $R \subset \dfr$ such that if $f \in R$ and $P$ is the set of
all hyperbolic periodic points of $f$, then both the sets $\cup_{p\in
P} W^s(p)$ and $\cup_{p\in P} W^u(p)$ are dense in $M$.

Furthermore, let $U \subset M$ be an open connected subset such that
it contains no periodic point for $f$ and suppose that, for some
hyperbolic periodic point $p$, $U \cap W^s(p) \neq \emptyset$, then
$W^s(p)$ is dense in $U$.
\label{thm1.1}\end{thm}

The main tool of our proof is the theory of prime ends \cite{Ca13}
\cite{Mather82a} and Lefschetz fixed point theorem. Also, Arnold's
conjecture for symplectic fixed points, as proved by Conley \& Zehnder
\cite{CZ83} for the case of torus, turned out to be essential for our
results on $\T^2$.

Our work was motivated in our effort to prove the so-called $C^r$
closing lemma, as conjectured by Poincar\'e \cite{Po1892}. Our result
provides a strong evidence supporting the $C^r$ closing lemma. We will
discuss this and other related problems in the end of the paper.

\vspace{1ex}
I am grateful to Kamlesh Parwani for useful comments and corrections
to the original manuscript.

\section{Prime ends and a Theorem of Mather}

Let $U$ be a simply connected domain of $S^2$ whose complement
contains more than one point. One can define the prime end
compactification of $U$, introduced by Caratheodory \cite{Ca13}, by
adding the circle $S^1$. Each point in $S^1$ are defined by a sequence
of nested arcs in $U$. $S^1$ inherits a topology from $U$. We refer to
Mather \cite{Mather82a} for definitions and general discussions. Here
we only summarize what we will need in this paper. The
compactification $\hat{U} = U \sqcup S^1$ of $U$, with inherited
topology on $S^1$, is homeomorphic to a closed disk $D$.  If $U
\subset S^2$ is an invariant subset for a homeomorphism on $S^2$, then
the function $f$ extends continuously to a homeomorphism $\hat{f}:
\hat{U} \ra \hat{U}$. Moreover, $\hat{f}$ restricted to the prime ends
$S^1$ is a circle homeomorphism.

When the rotation number $\rho$ of $\hat{f}|_{S^1}$ is rational, then
$\hat{f}$ has a periodic points on the prime ends. But it may happen
in general that $f$ does not have any periodic points on $\partial
U$. However, if $f$ is area-preserving, then this will never happen.

\begin{lem}
Let $f$ be an area-preserving homeomorphism on $S^2$ and let $U$ be a
simply connected, invariant set whose complement contains more than
one point. Let $\hat{f}: \hat{U} \ra \hat{U}$ be the extension of $f$
to the prime end compactification of $U$. If $\hat{f}$ has a periodic
point in its prime ends then $f$ has a periodic point in $\partial U$.
\label{lempe}\end{lem}

The proof is simple, we refer to Mather \cite{Mather82a} and Franks \&
Le Calvez \cite{FL03}. The next theorem, which is due to Mather
\cite{Mather82a} states that in generic situations, there are no
periodic points on the prime ends.

\begin{thm}
Let $f$ be an area-preserving diffeomorphism on $S^2$ and let $U$ be a
simply connected, invariant set whose complement contains more than
one point. Further assume that $f$ satisfies the generic conditions
G1, G2 and G3. Then there is no periodic point for $f$ on the boundary
of $U$ and as a consequence, there is no periodic point for $\hat{f}$
on the prime ends and the rotation number of $\hat{f}$ on the prime
ends is irrational.
\label{thmpe}\end{thm}

The proof is based on the following ideas: by generic conditions G1
and G2, $f$ can not have any elliptic periodic points on $\partial
U$. If there is any hyperbolic periodic points on $\partial U$, then
the boundary of $U$ has to contain some branches of the stable and
unstable manifolds of the hyperbolic periodic points and this will
lead to a contradiction to the generic property G3. For details of the
proof, see Mather \cite{Mather82a} and Franks \& Le Calvez
\cite{FL03}.

The theory of prime ends can be easily generalized to other orientable
compact surfaces $M$, possibly with boundaries. Let $U \subset M$ be a
connected open domain in $M$. Assume that the boundary of $U$ contains
finite number of connected pieces and each piece contains more than
one point. Then $U$ can be compactified by adding prime ends. In this
case, the prime ends are homeomorphic to a union of finitely many
circles $S^1$. The number of circles $S^1$ added is the number of
boundary pieces seeing from the inside of $U$. It may be more than the
number of connected components of $\partial U$.  The compactification
of $U$, denoted as $\hat{U}$, is a compact surface with boundaries,
each $S^1$ will be a hole for $\hat{U}$. If $f: M \ra M$ is a
homeomorphism and $U$ is invariant for $f$, then $f$ extends
continuously to a homeomorphism on $\hat{U}$. We denote this extension
by $\hat{f}$.

Lemma \ref{lempe} and Theorem \ref{thmpe}, with obvious modifications,
are both true for $\hat{U}$.

In case where $U$ has infinitely many connected boundary components,
or where some of the boundary components are single points, one can do
a prime end compactification to finitely many isolated connected
boundary pieces that contain more than one point. Again, the number of
boundary pieces are defined to be the number seeing from inside $U$

\section{Periodic-point-free regions}

Let $R$ be a residual subset of $\dfr$ such that for any map in $R$ the
generic properties G1, G2, G3 and G4 hold. Whenever G4 is assumed, we
require the map to have sufficient smoothness, say $r \geq 16$. Our
final result will be true for any $r=1, 2, \ldots,$.

Fix a map $f \in R$. Let $U \subset M$ be a connected open set such that
$f$ has no periodic point in $U$. Our goal is to show that either $U$ is
filled by stable and unstable manifolds of some hyperbolic periodic points,
or by an arbitrarily small $C^r$ perturbation of $f$, we can create a
periodic point in $U$.

Let $A_U$ be the set of all points whose orbit under $f$ passes
through $U$. i.e., $$A_U = \cup_{i= -\infty}^\infty f^i(U).$$ Since
$f$ is area-preserving, almost every point is recurrent. Since $M$
has finite area, $A_U$ has finite number of connected components. Each
of these components are open and periodic. Without loss of generality,
we may assume that $A_U$ has only one component and this component
certainly contains $U$. The case with more than one components can be
reduced to this case simply by considering powers of $f$.

The set $A_U$ is invariant under $f$ and it is periodic point free.
We are interested in whether the closure of $A_U$ contains any
periodic point. The set $A_U$ typically has infinitely many connected
boundary components, so we can not apply Lemma \ref{lempe} and Theorem
\ref{thmpe}, or other results obtained by the prime ends directly. We
have the following simple lemma.

\begin{lem}
  The closure of $A_U$, $\overline{A_U}$, is a compact invariant set which
  contains no elliptic periodic point for $f$. Here we assume $f$ has
  sufficient smoothness so that the generic condition G2 holds.
\label{lem3.1}\end{lem}

\noindent {\em Proof}: Recall that the map $f \in R$ satisfies the
generic property G2, thus all elliptic periodic points are Moser
stable. This implies that surrounding each elliptic periodic point
there are infinitely many invariant curves. If $\overline{A_U}$
contains an elliptic periodic point, then $A_U$ has to intersect an
invariant curve with irrational rotation number close to the elliptic
point. The set $A_U$ is invariant, thus it completely covers that
curve. By Birkhoff fixed point theorem, or more generally the
Aubry-Mather theory, any neighborhood of the invariant curve contains
infinitely many periodic orbits. This contradicts to the assumption
that $A_U$ is periodic point free.

This proves the lemma.

The closure of $A_U$ may still contain hyperbolic periodic points. We
show that this happens only if $A_U$ contains stable or unstable
manifolds of hyperbolic periodic points.

\begin{lem}
If the closure of $A_U$, $\overline{A_U}$, contains a hyperbolic periodic
point $p$, then $U \cap W^s(p) \neq \emptyset$ and $U \cap W^u(p) \neq
\emptyset$.
\label{lem3.2}\end{lem}

\noindent {\em Proof}: Recall that for $f \in R$, generic property G4
implies that $\overline{W^s(p)} = \overline{W^u(p)}$, it suffices to
show that either $A_U \cap W^s(p) \neq \emptyset$ or $A_U \cap W^u(p)
\neq \emptyset$.

First, we remark that if $M=S^2$ or $M=\T^2$ and generic property G5
is assumed, then the lemma can be easily. In this case, each branch of
the stable and unstable manifolds of $p$ intersect and the
intersection is transversally, there is a homoclinic tangle in the
neighborhood of $p$ and this homoclinic tangle divides a small
neighborhood of $p$ into infinitely many small rectangles. Choose a
small neighborhood $p$ such that all of these rectangles have areas
smaller than that of $U$. If $p$ is in the closure of $A_U$, then, for
some integer $i$, $f^i(U)$ intersects one of these small
rectangles. By the area preserving property of $f$, $f^i(U)$ can not
be totally contained in any single rectangle and thus it has to
intersect the boundary of the rectangle. This implies that $f^i(U)$
contains a point on either the stable manifold or the unstable
manifold of $p$.

For the general case, suppose that there is a hyperbolic periodic
point $p \in M$ such that $A_U \cap W^s(p) = \emptyset$, we want to
show that $p$ is not in the closure of $A_U$. Let $V$ be the connected
component of $M \backslash \overline{W^s(p)}$ containing $A_U$. Since
the complement of $V$ is connected, the boundary of $V$ consists of
finitely many connected pieces and each piece contains more than one
point. Therefore $V$ has a prime end extension, $\hat{V}$ with the
extended map $\hat{f}$ on $\hat{V}$. By Mather's theorem
\cite{Mather82a}, which is proved in this general setting, $\hat{V}$
has no periodic points on the boundary of $\hat{V}$. This implies that
$p \notin \overline{V}$ and since $A_U \subset V$, this implies that
the closure of $A_U$ contains no periodic point.

This proves the lemma.

To get rid of isolated points in the boundary of $A_U$, we let
$A=\mbox{int} (\bar{A_U})$. Then any isolated connected boundary piece
contains more than one point. If the closure of $A_U$ is periodic
point free, then so is the closure of $A$.

\section{Isotopic to identity cases}

In this section we prove our main lemma in the case where the map is
isotopic to identity.

\begin{lem}
  Let $f \in R$ be a generic diffeomorphism of $M$ and $f$ is isotopic
  to identity. Let $A \subset M$ be an open, connected, periodic
  point free, $f$-invariant set whose closure contains no periodic
  point. We further assume that $A = \mbox{int} \bar{A}$. Then either
  $A= M= \T^2$ or $A$ is homeomorphic to an open annulus and its prime
  end extension is a closed annulus. \label{mainlemma}
\end{lem}

\noindent {\em Proof}: We first assume that $M \neq \T^2$. Since the
Euler characteristic of $M$ is non-zero, by Lefschetz fixed point
theorem, $f$ contains at least one fixed point. This implies that $A
\neq M$. Let $\mbox{Fix}(f)$ be the set of all fixed point for $f$. It
contains finite number of points and it is contained in the interior
of $M \backslash A$. Let $B \subset M$ be the union of the connected
components of $M\backslash A$ intersecting $\mbox{Fix(f)}$. The set
$B$ is closed and has finite number of connected components and each
contains at least one fixed point of $f$. Let $C = M \backslash B$,
obviously $C$ is open, $A \subset C$ and the closure of $C$ contains
no fixed point. Let $\hat{C}$ be the prime end extension of $C$ and
let $\hat{f}: \hat{C} \ra \hat{C}$ be the extension of $f$. Then
$\hat{f}$ has no fixed point. The set $\hat{C}$ is a compact surface
with boundary, its topology is uniquely determined by the number of
handles and the number of holes on $S^2$. Let the number of handles of
$\hat{C}$ be $k$ and the number of holes be $l$. The map $\hat{f}$
keeps invariant of all the boundary pieces. The Euler characteristic
number of $\hat{C}$ is $2-2k-l$. Even though $\hat{f}$ may not
neccessarily be homotopic to identity, but its induced map on homology
is identity. Hence its Lefschetz number is the same as its Euler
characteristic number. Since $\hat{f}$ has no fixed point, this
implies that $2-2k-l=0$. This implies that either $k=0$, $l=2$ or
$k=1$, $l=0$. The latter implies that $C$ is a torus and thus $M$ is a
torus, which we assumed was not the case. The first case implies that
$\hat{C}$ is a closed annulus.

We claim that $A=C$. If not, there is a point $x \in C\backslash
\bar{A}$. i.e., $x$ is in the interior of $C \backslash A$. For if no
such point exists, then $C \subset \bar{A}$. But $C$ is open, this
implies that $C \subset \mbox{int}\bar{A} = A$, by the assumption on
$A$. Let $B_x$ be the connected component in $C \backslash \bar{A}$
containing $x$. Then $B_x$ is closed and $B_x \cap \partial C =
\emptyset$. For if $B_x \cap \partial C \neq \emptyset$, then $B_x
\subset B$, where $B$ is the set of connected components containing
fixed points in the complement of $A$, as defined previously. Since
$B_x$ has positive area, it must be periodic under $f$. Let $k$ be the
period of $B_x$. Let $\mbox{Fix}(f^k)$ be the set of all periodic
points of $f$ with period $k$. Again, $\mbox{Fix}(f^k)$ contains only
finite number of points. Let $B_k$ be the union of all connected
components of $C \backslash A$. The set $B_k$ has finite number of,
say $l$, connected components and all components are periodic with
period $k$. Moreover, $B_k \cap \partial C = \emptyset$. Let $C_x = C
\backslash B_k$, then $C_x$ is an open set and $C_x =
\mbox{int}\bar{C}_x$. The prime end extension of $C_x$ is a closed
annulus with $l$ interior disks removed. and there is no periodic
points of period $k$ in $C_x$. However, the Lefschetz number for $f^k$
on $C_x$ is $l \neq 0$, this implies that $f^k$ has to have at least
one fixed point in $C_x$. This contradiction proves that $A=C$. i.e.,
$A$ is an annulus.

If $M=\T^2$, then either $A=T^2$ or there is a point $x \in
\mbox{int}(M \backslash A)$. Let $B'$ be the connected components
containing $x$ in the complement of $A$. Again, $B'$ must be periodic
and let this period be $k$. Let $B$ be the union of all the connected
components of $M \backslash A$ that are fixed under $f^k$. $ \emptyset
\neq B' \subset B$.  Let $C = A \backslash B$. Then this case follows
from above arguments by considering $f^k$ on $C$.

This proves our lemma.

\section{A lemma on periodic points}

Now we consider the case where $f$ is not necessarily isotopic to
identity. To prove the same result, we need to prove a simple theorem
on the existence of periodic points for maps of compact surfaces with
boundaries.

\begin{lem}
Let $M$ be a compact, connected surface, possibly with
boundary. Assume that $M$ has nonzero Euler characteristic. Suppose
that $f: M \rightarrow M$ is a homeomorphism, then $f$ has a periodic
point. Moreover, for any positive integer $n$, there exists infinitely
many positive integers $i$ such that the Lefschetz number $L(f^{ni})$
is negative.
\label{perlem}
\end{lem}

\noindent {\em Proof}: Let $L(f)$ be the Lefschetz number of the map
$f$. By Lefschetz fixed point theorem, it suffices to show that
$L(f^n)$ is nonzero for some nonzero integer $n$.

We consider homology groups of $M$ with real coefficients. If $H_2(M)
\neq 0$, by taking an interate of $f$, we may assume that the induced
map on $H_2(M)$ is the identity. Let $tr_1(f_*)$ be the trace of
the induced map $f_*: H_1(M) \rightarrow H_1(M)$, then $L(f^n) = 2-
tr_1(f^n_*)$ if $M$ is orientable and without boundary, and $L(f^n) =
1- tr_1(f^n_*)$ if $M$ is non-orientable or has boundary. The theorem
follows easily for the cases where the dimension of $H_1(M)$ is zero
(the sphere and the projective plane). The torus, the Klein bottle,
the M\"obius strip and the annulus all have Euler characteristic zero
and hence are excluded. The only cases we need to consider are where
the dimension of the first homology $H_(M)$ are greater than or equal
to 3. Let $\lambda_1, \lambda_2, \ldots, \lambda_l$ be the eigenvalues
of the induced isomorphism $f_*: H_1(M) \rightarrow H_1(M)$, where $l$
is the dimension of $H_1(M)$. We can write $\lambda_i= r_i \alpha_i$,
where $r_i$ is a positive real number and $\alpha_i$ is complex number
on the unit circle, for $i=1, 2, \ldots, l$. There exists a sequence
of integers $\{n_k\}_{k=1}^\8$ such that $\alpha_i^{n_k} \ra 1$ as $k\ra
\8$, for all $i=1, 2, \ldots, l$.

If $r_i > 1$ for some $i=1, 2, \ldots, l$, then $tr_1(f^{n_k}_*) =
\sum_{i=1}^l \lambda_i^{n_k} \rightarrow \8$ as $k \ra \8$. This
implies that for large $k$, $L(f^{n_k}) \neq 0$, the lemma follows. If
$r_i < 1$ for some $i=1, 2, \ldots, l$, then we consider $f^{-1}$, the
inverse of $f$. The same argument shows that $L(f^n) \neq 0$ for some
negative integer $n$. Since  $L(f^n) = L(f^{-n})$, the lemma again
follows. Finally, if $r_i =1$ for all $i = 1, 2, \ldots, l$, then
$tr_1(f^{n_k}_*) = \sum_{i=1}^l \lambda_i^{n_k} \rightarrow l$ as $k
\ra \8$. As $j \geq 3$, $2- tr_1(f^{n_k}_*) \neq 0$ and $1 -
tr_1(f^{n_k}_*) \neq 0$ for large $k$. This shows that for infinitely
many positive integer $k$, the Lefschetz number $L(f^{n_k})$ is
negative. Lefschetz fixed point theorem concludes that there are at
least one fixed point for $f^{n_k}$ for such $k$.

For any positive integer $n$, replacing $f$ with $f^n$, the above
arguments show that there are infinitely many positve integer
$i$ such that $L(f^{ni})$ is negative.

This proves the lemma.

For orientable compact surfaces, only $\T^2$ and annulus have zero
Euler characteristics.  Lemma \ref{mainlemma} would follow from the
above lemma if the periodic point free set $A$ in Lemma
\ref{mainlemma} has a prime end extension that makes it into a compact
manifold with boundary. We have to show that $A$ can't have infinitely
many boundary pieces. The proof of the lemma is basically a
verification of this fact.

\section{General case}

Let $f \in R$ be a generic diffeomorphism on the compact surface
$M$. Let $A \subset M$ be an open, connected, periodic point free,
$f$-invariant set whose closure contains no periodic point. Assume $A
= \mbox{int}(\bar{A})$. For any integer $i$, we let $\mbox{Fix}(f^i)$
be the set of all periodic points of $f$ with period $i$. If $A \neq
M$, $M\backslash A$ is nonempty and contains at least one connected
component. The number of the connected components in $M \backslash A$
may be infinite.  Let $B_i$ be the set of all connected components of
$M \backslash A$ containing a periodic point of period $i$. $B_i$ is
closed and $B_i \cap \mbox{Fix}(f^i) \neq \emptyset$, if
$\mbox{Fix}(f^i) \neq \emptyset$. $B_i$ has finite number of connected
components. For any positive integers $i, j$, each component of $B_j$
is either a component of $B_i$ or disjoint from $B_i$. If $j =ki$ for
some positive integer $k$, then $B_i \subset B_j$.

Let $C_i = M \backslash B_i$. Then $C_i$ has finite number of boundary
components. Let $\hat{C}_i$ be the prime end extension of $C_i$ and it
is a compact surface with boundary. Let the number of handles of
$\hat{C}_i$ be $m$ and the number of holes be $n$. For any integer
$k$, the number of handles of $\hat{C}_{ki}$ is smaller than or equal
to $m$. Therefore, there exists an integer $i^*$ such that the number
of handles of $\hat{C}_{ki^*}$ is a constant for all positive integer
$k$.

If $\mbox{Fix}(f^i)$ is empty for all positive integer $i$, then
either $A=M$, then we set $C_i=M$ for all $i$, or $A \neq M$, then we
pick a point $x \in M\backslash A$ and let $B_i$ be the connected
component of $M \backslash A$ and its iterates under $f$. Since $f$ is
area preserving, $B_i$ has finite number of components.

By our construction, $\hat{f}^i: \hat{C}_i \ra \hat{C}_i$ has no
fixed points for all $i \geq 1$.

We are ready to prove  Lemma \ref{mainlemma} for the general case.

\begin{lem}
  Let $f \in R$ be a generic diffeomorphism of $M$. Let $A \subset M$
  be an open, connected, periodic point free, $f$-invariant set whose
  closure contains no periodic point. We further assume that $A =
  \mbox{int} \bar{A}$. Then either $A= M= \T^2$ or $A$ is homeomorphic
  to an open annulus and its prime end extension is a closed annulus.
\label{mainlemmag}\end{lem}

\noindent {\em Proof}: First if there exists a positive integer $k$
such that $f^k$ is isotopic to identity, the lemma is proved in the
same way as Lemma \ref{mainlemma}. One just need to consider $f^k$.

We first suppose that $M \neq \T^2$.

As described before, there exists an integer $i^*$ such that the
number of handles of $\hat{C}_{ki^*}$ is a constant for all positive
integer $k$. We now consider $\hat{f}^{i^*}: \hat{C}_{i^*} \ra
\hat{C}_{i^*}$. It is a homeomorphism on compact surface, keeping
invariant of all boundary pieces. Suppose $\hat{f}^{i^*}:
\hat{C}_{i^*} \ra \hat{C}_{i^*}$ has no periodic points, then by Lemma
\ref{perlem}, $C_{i^*}$ is an open annulus and $\hat{C}_{i^*}$ is a
closed annulus. We claim that $A = C_{i^*}$. Suppose not, then there
is a point $x \in C_{i^*}$, but $x \notin A$. Since $A = \mbox{int}
\bar{A}$ and by the definition of $C_{i^*}$, $x \notin \partial
A$. Let $C_x$ be the connected component of $C_{i^*}\backslash
\bar{A}$ containing $x$. Then $C_x$ an open disk in a annulus. Since
$f$ is area preserving, $C_x$ must be periodic and it contains a
periodic point. This contradicts to the assumption that $C_{i^*}$ has
no periodic point.

Now suppose that $\hat{f}^{i^*}: \hat{C}_{i^*} \ra \hat{C}_{i^*}$ has
a periodic point $x$ with period $p$. By the assumption on $A$, $x
\notin \bar{A}$. By Lemma \ref{perlem}, there is a positive integer
$j$ such that the Lefschetz number of the $pj$'s iterate of
$\hat{f}^{i^*}$ on $C_{i^*}$ is negative. i.e.,
$L((\hat{f}^{i^*})^{pj}| C_{i^*}) < 0$, where $p$ is the period of
$x$. Now consider the prime end extension of $C_{i^*pj}$. By our
choice of $i^*$, $\hat{C}_{i^*pj}$ is topologically $\hat{C}_{i^*}$
with finitely many, say $k$, $k >0$, open disks removed and these
disks are periodic with period $pj$ under the extended map
$\hat{f}^{i^*}$. We have the following relations on the Lefschetz
numbers $L(\hat{f}^{i^*pj}| C_{i^*pj}) = L((\hat{f}^{i^*})^{pj}|
C_{i^*}) - k <0$. Lefschetz fixed point theorem implies that
$\hat{f}^{i^*pj}$ has a fixed point on $C_{i^*pj}$, which is
impossible by the definition of $C_{i^*pj}$. This contradiction shows
that $\hat{f}^{i^*}: \hat{C}_{i^*} \ra \hat{C}_{i^*}$ has no periodic
point and hence $A=\hat{C}_{i^*}$ is an annulus.

We are left with one case: $M=\T^2$ and $A \neq M$. The same argument
works in this case too.

This proves the lemma

\section{Maps on annulus}

Let $U$ be a connected open subset of $M$ and let $A_U = \cup_{i =
  -\infty}^\infty f^i(U)$. If the closure of $A_U$ contains no
periodic point, then $A=\mbox{int}(\bar{A}_U)$ is open, containing no
periodic points in its closure. By the above lemma, if $M \neq \T^2$,
then $A$ is a union of finite disjoint open annuli, periodic under
$f$. The dynamics on Annulus have been well studied (cf.\ Franks
\cite{Franks96}, Le Calvez \& Yoccoz \cite{LY97} \ The following lemma
shows that, if $A$ is an annulus, we can perturb $f$ with an
arbitrarily small $C^r$ perturbation to create a periodic point in
$U$. The same result was also used in \cite{FL03}

\begin{lem}
Fix $f\in R$ and assume $M \neq \T^2$. Let $U$ be a connected open
subset of $S^2$ and $A_U = \cup_{i = -\infty}^\infty f^i(U)$. Assume
that $U$ does not intersect any stable or unstable manifolds of
hyperbolic periodic points of $f$.  Then for any $C^r$ neighborhood
$V$ of $f$, there exists $g \in V$ such that the support of $g-f$ is
contained in the interior of the closure of $A_U$ and $g$ has a
periodic point in $U$.
\label{lemannulus}\end{lem}

\noindent {\em Proof}: Since $U$ does not intersect any stable or
unstable manifolds of hyperbolic periodic points of $f$, By above
lemmas, $f$ has no periodic point in $A_U$ and no periodic point in
the closure of $A_U$ and therefore, $A= \mbox{int} (\bar{A}_U)$ is a union
of finite disjoint periodic annuli. Without loss of generality, we may
assume that $A_U$ itself is an open annulus.

Using prime end extension, we obtain an area-preserving continuous map on
the prime end closure of $A_U$, still denoted by $\overline{A_U}$. Let $(x,
y), x \in \mathbb{R} \; (\mod 1), y \in [0, 1]$ be a coordinate on
$\overline{A_U}$.  Since $f$ preserves invariant measure $\mu$, by Birkhoff
Ergodic Theorem, for $\mu$-almost every point $z=(x, y) \in
\overline{A_U}$, the rotation number $$\rho(z, f) = \lim_{i \rightarrow
\infty} \frac{\pi_x \tilde{f}^i(z)}{i}$$ is well defined. Here $\pi_x$ is the
projection on $\overline{A_U}$ into its first coordinate and $\tilde{f}$ is
a lift of $f$ to its universal cover $\mathbb{R} \times [0, 1]$. A different
lift of $f$ yields a different rotation number that differs by an integer.

Since there is no periodic points in $\overline{A_U}$, by Franks
theorem \cite{Franks96}, there exists an irrational number $\alpha \in
\mathbb{R}$ such that for almost all $z \in A_U$, the rotation number
$\rho(z, f)$ exist and $\rho(z, f) =\alpha$. In particular, if $z$ is
in the boundaries of $A_U$ then $\rho(z, f) =\alpha$.

Let $\gamma$ be a simple closed curve in the interior of $A_U$. We may
assume that $\gamma$ is homotopically non-trivial in $A_U$. Take a small
tubular neighborhood $\gamma_\delta$ of $\gamma$ in the interior of $A_U$
and parametrize this tubular neighborhood by $\gamma_\delta: S^1 \times
[-\delta, \delta] \rightarrow: A_U$ for some small $\delta$. In fact, for
convenience we may even assume that $\gamma_\delta$ is area-preserving. Let
$\beta: [-\delta, \delta] \rightarrow \mathbb{R}$ be a $C^\infty$ function
such that $\beta(t) >0$ for all $-\delta < t < \delta$ and $\beta(-\delta) =
\beta(\delta) =0$ and $\beta$ is $C^\infty$ flat at $\pm \delta$. i.e., all
the derivatives of $\beta(t)$ at $\pm \delta$ are zero.

Let $h_\ep: A_U \rightarrow A_U$ be a $C^\infty$ diffeomorphism such
that if $z \notin \gamma_\delta$, $h_\ep(z) =z$ and if $z \in
\gamma_\delta$, $h_\ep = \gamma_\delta \circ T_\ep \circ
(\gamma_\delta)^{-1}$ where $T_\ep ( \theta, t) = (\theta + \ep
\beta(t), t)$ for all $\theta \in S^1$ and $t \in [-\delta,
\delta]$. We remark that $h_\ep \rightarrow \mbox{Id}$ in $C^\infty$
topology as $\ep \rightarrow 0$ and the mean rotation number for
$h_\ep$ with respect to the area $\mu$ is
$$\rho_\mu (A_U, h_\ep) = \frac{1}{\mu(A_U)} \int_{-\delta}^\delta \ep
\beta(t) dt$$

Therefore for any $\ep > 0$, $\rho_\mu(A_U, h_\ep \circ f) = \rho_\mu
(A_U, h_\ep) + \rho_\mu(A_U, f) > \alpha$, this implies that there
exists a point $y^* \in A_U$, such that $\rho(y^*, h_\ep \circ f) >
\alpha$. Since $\rho(z, h_\ep \circ f) = \rho(z, f) = \alpha$ for all
$z \in \partial \overline{A_U}$, we conclude, from Franks' theorem
\cite{Franks96}, that for any rational number $p/q$, such that $\alpha
< p/q < \rho(y^*, h_\ep \circ f)$, there exists a periodic point of
period $p/q$ for the map $h_\ep \circ f$.

Thus, for any $\ep > 0$, there are infinitely many periodic points for
$h_\ep \circ f$ in the interior of $A_U$. In fact, all of these periodic
points have to pass through the strip $\gamma_\delta$. However, these
periodic points may be far away from $U$. To find periodic points in $U$,
we need to do some estimates on these orbits.

Since $A_U = \cup_{i \in \mathbb{Z}} f^i(U)$, for any point $z
\in \gamma_\delta \subset A_U$, there exists an integer $n_z \in
\mathbb{Z}$ and a neighborhood of $z$, $W_z \subset A_U$, such
that $f^{n_z} (W_z) \subset U$. $\{W_z, \; z \in \gamma_\delta\}$ forms an
open cover for the compact set $\gamma_\delta$. Let $W_{z_1}, W_{z_2},
\ldots, W_{z_k}$ be a finite subcover of $\gamma_\delta$ and let $N= \max
\{|n_{z_1}|, |n_{z_2}|, \ldots, |n_{z_k}| \}$.

The integer $N$ chosen above has a very important property: for any $z \in
\gamma_\delta$, the orbit segment $\{ f^{-N}(z), f^{-N+1}(z), \ldots,
f^N(z)\}$ intersects $U$ at least once. Or equivalently, the set
$\cup_{i=-N}^N f^i(U)$ covers $\gamma_\delta$. Since $U$ is open, this same
property holds for all $g$ sufficiently close to $f$ in $C^0$
topology. i.e., the orbit segment $\{ g^{-N}(z), g^{-N+1}(z), \ldots,
g^N(z)\}$ intersects $U$ for all $z \in \gamma_\delta$, provided that $g$
is sufficiently close to $f$.

Above arguments show that if $\ep >0$ is small enough, $h_\ep \circ f$ has
infinitely many periodic orbits and all of these periodic orbits intersect
$U$. 

This proves the lemma.

\section{Maps on torus and Arnold's conjecture}

The final case is where $M=\T^2$ and $f$ has no periodic point. We
will show that such $f$ is not generic and it can be perturbed to
create a periodic point.

Let $f_{*1}: H_1(\T^2, \R) = \R^2 \ra \R^2$ be the induced map on the
first homology of $\T^2$. Let $\lambda_1$, $\lambda_2$ be the
eigenvalues of $f_{*1}$, $\lambda_1 = \bar{\lambda}_2$. Since $f$ has
no periodic point, the Lefschetz number $L(f^k) =0$ for all $k$. This
implies that $\lambda_1^k + \lambda_2^k =2$, for all $k$. We must have
$\lambda_1 = \lambda_2 =1$. Since $f_{*1}$ and its inverse are both
integer matrices, we have only two choices: $f$ is isotopic to
identity, where $f_{*1} = I$  or $f$ is isotopic to a Dehn twist,
i.e., for some integer $k \neq 0$,
$$f_{*1} = \begin{pmatrix}1 & k \cr 0 &1 \end{pmatrix}$$

We first consider the case where $f$ is isotopic to identity. The
proof is basically an application of the Arnold conjecture
\cite{Arnold78} as proved by Conley and Zehnder \cite{CZ83}.

\begin{lem}
Let $f: \T^2 \ra \T^2$ be an area-preserving diffeomorphism such that
$f$ is isotopic to identity. Then for any $C^r$ neighborhood $V$ of
$f$, there exists $g \in V$ such that $g: \T^2 \ra \T^2$ has a
periodic point.
\label{lemtorus}\end{lem}

\noindent {\em Proof}: 

Let $\tilde{f}: \R^2 \ra \R^2$ be a lift of the map $f$ on
$\T^2=\R^2/\Z^2$ to its universal cover $\R^2$. Let $\pi_i$, $i=1, 2$
be the projection of $\R^2$ to its first and second coordinates
respectively. The average rotation numbers for the map $\tilde{f}$,
are defined to be
$$\rho_i(\tilde{f}) = \int (\pi_i(\tilde{f}(p)) - \pi_i(p)) d\mu.$$
The combination $\rho(f) = (\rho_1(\tilde{f}), \rho_2(\tilde{f})),
\mod \Z^2$ is called the average (or mean) rotation vector for
$f$. The rotation vectors are well defined for maps isotopic to
identity.

We want to do a small perturbation to $f$ so that each component of
the average rotation vector is a rational number. This is easy: one
composes the map $f$ with $T_{(\ep_1, \ep_2)}(x, y) = (x_1 +
\ep_1, x_2 +\ep_2)$, then $\rho(f\cdot T_{(\ep_1, \ep_2)}) = \rho(f) +
(\ep_1, \ep_2), \mod \Z^2$. By properly choosing small $\ep_1$ and
$\ep_2$, we obtain a rational rotation vector for $f\cdot T_{(\ep_1,
  \ep_2)}$.

There is a positive number $i$ such that the mean rotation vector for
$( f\cdot T_{(\ep_1, \ep_2)})^{ik}$ is an integer vector, which is
equivalent to zero on the torus. $( f\cdot T_{(\ep_1, \ep_2)})^{ik}$
is isotopic to identity and preserves a smooth area element. The
Arnold conjecture, as proved by Conley and Zehnder \cite{CZ83} in the
case of torus, implies that $( f\cdot T_{(\ep_1, \ep_2)})^{ik}$ has at
least three fixed points, four if all non-degenerate. This implies
that $f\cdot T_{(\ep_1, \ep_2)}$ has periodic points of period $ki$.

This proves the lemma.

We use Poincar\'e-Birkhoff Theorem for the case where $f$ is isotopic
to a Dehn twist.

\begin{lem}
Let $f: \T^2 \ra \T^2$ be an area-preserving diffeomorphism such that
$f$ is isotopic to a Dehn twist. Then for any $C^r$ neighborhood $V$
of $f$, there exists $g \in V$ such that $g: \T^2 \ra \T^2$ has
infinitely many periodic points.
\end{lem}

Let $((x, y), \mod \Z^2)$ be a coordinate systems on $\T^2$ such that
$f$ is isotopic to the map $(x, y) \mapsto ((x +ky, y), \mod \Z^2)$
for some non-zero integer $k$. We may assume, without loss of
generality, that $k >0$. In this homotopy class, there are maps
without any periodic point. For example, the map $(x, y) \mapsto (x
+ky, y+ \alpha)$ has no periodic point if $\alpha$ is
irrational. However, there are periodic points when $\alpha$ is
rational. Our first step is to perturb the map so that it has a
rational vertical rotation number.

Lift the map $f$ in the $y$ direction, we obtain a map on the infinite
cylinder $\tilde{f}^y: S^1 \times \R \ra S^1 \times \R$. Define the
mean vertical rotation number
$$\rho_2(\tilde{f}^y) = \int (\pi_2(\tilde{f}^y(p)) - \pi_2(p))
d\mu,$$ where $\mu$ is the area element. $\pi_2(\tilde{f}^y(p)) -
\pi_2(p)$ is independent of choices of the covering points. We define
the mean vertical rotation number of $f$ to be $\rho_2(f) =
\rho_2(\tilde{f}^y), \mod 1$.  This is independent of the lift.

By composing $f$ with the map $(x, y) \mapsto (x, y + \ep)$ for some
small $\ep$, we obtain a map $g$ on the torus such that its mean
vertical rotation number is rational. This implies that there is a
positive integer $l$ such that the vertical rotation number of $g^l$
is zero. Now, choose a lift of $g^l$, $G^y$, in the $y$ direction such
that its mean vertical rotation number is zero (instead of being a
non-zero integer). Then $G^y$ is an area preserving map on
the infinite cylynder $S^1 \times \R$ which is also exact, i.e.,
integral of the 1-form $ydx$ over any closed curve on the cylinder is
invariant under the map. Moreover, if we let $G$ be the lift
of $G^y$ to $\R^2$, we have that $\pi_1(G(p)) - x
- kl y$ is uniformly bounded for any $p = (x, y) \in \R^2$.

Poincar\'e-Birkhoff twist map theorem implies that there are at least
two fixed points and infinitely many periodic points for
$G^y$. This implies that there are infinitely many periodic
points for $g$ on $\T^2$. 

This proves the lemma.

The Poincar\'e-Birkhoff twist map theorem (or Poincar\'e's last
geometric theorem) shows the existence of fixed points for
area-preserving maps on the annulus with twist condition. Here we have
an infinite cylinder with infinite twists on two ends. There are two
ways to work around this. One is to directly apply Birkhoff's proof
\cite{Birkhoff66} to the annulus $\{|y| \leq M\}$ with large $M >0$. The
annulus is not invariant but the proof works in the same way, as long
as one has exactness in the area preserving property. Another way to
prove the result is to modify the map so that all horizontal lines are
fixed for large values of $|y|$. Again, the exactness is
neccessary. Since these techniques are well known, we will not give
details here.

\section{Proof of the main theorem}

Let $M$ be a compact surface and let $R \subset \dfr$ be the set of
area-preserving diffeomorphisms on $M$ satisfying G1-G4. We first
assume that $r \geq 16$. In addition, we assume that $R$ satisfies the
following two conditions: for any $f \in R$,

\begin{itemize}
\item[G6] every invariant open annulus contains a
periodic point;

\item[G7] if $M= \T^2$, then there is a periodic point.
\end{itemize}

By Lemma \ref{lemtorus}, G7 is an open and dense condition, hence
generic. By Lemma \ref{lemannulus}, G6 is a generic
condition. Therefore, $R$ is a residual set. We claim that for any $f
\in R$, the stable and unstable manifolds of hyperbolic periodic
points are dense in $M$. Suppose this is not true and there exists an
open set $U \subset M$ such that $U$ does not intersect stable
manifold and unstable manifold of any hyperbolic periodic point. Then
neither does the invariant set $A_U=\cup_{i \in \Z} f^i(U)$. Then by
Lemma \ref{lem3.1} and Lemma \ref{lem3.2}, the closure of $A_U$
contains no periodic points. Let $A = \mbox{int} (\overline{A_U})$,
then the closure of $A$ contains no periodic point. Lemma
\ref{mainlemmag} shows that $A$ must be an open annulus or $M$ must be
a torus. This contradicts to conditions G6 and G7.

This proves the first part of our theorem for $r \geq 16$. For lower
smoothness, i.e., for $r = 1, 2, \ldots, 15$, we first note that the
residual set $R \subset \dfr$, $r \geq 16$ constructed above is dense
in $\dfr$ with $ r= 1, 2, \ldots, 15$. Moreover, for each open set $U
\subset M$, if $U$ intersects a piece of stable (or unstable) manifold
of a hyperbolic fixed point for some $f \in \dfr$, then same is true
for any map close to $f$. Since $M$ has a countable basis of open
sets, there is a residual subset $R^r \subset \dfr$ for all $r=1, 2,
\ldots, $ such that every open set intersects a piece of stable (and
unstable) manifold of some hyperbolic periodic point.

This proves the first part of our main theorem. For the second part of
the theorem, we need the following lemma.

\begin{lem}
Every prime end is an accumulation point of periodic points for a
generic $C^r$ surface diffeomorphism. More precisely, there is a
residual subset $R \subset \dfr$ (This set can be chosen to be the
same as above) such that for any $f \in R$, we have the following
property: Let $V$ be an open, $f$ invariant set $V$ with finite number
of connected boundary pieces and each boundary piece containing more
than one point, let $\hat{V}$ be the prime end extension of $U$, let
$z \in \hat{V}$ be a prime end, then there is a sequence of periodic
points of $f$, $\{p_n\}_{i=1}^\8$ such that $p_n \rightarrow z$ as $n
\rightarrow \8$.
\end{lem}

We will not give a detailed proof of this lemma. The proof follows
from Corollary 8.9 in Franks \& Le Calvez \cite{FL03}, which uses
Conley index in a small neigborhood of the prime end circle to obtain
periodic points (cf.\ Franks \cite{Franks99} and Le Calvez \& Yoccoz
\cite{LY97}). We remark that even though their results are on $S^2$,
there is no difference nearby one piece of prime ends. We also remark
that each prime end circle has irrational rotation, if a sequence of
periodic points approaches the prme end circle, then the sequence of
periodic orbits approaches every point in the prime end circle.

Now let $f \in R$ and $U$ be an open connected that contains no
periodic point for $f$. Assume that $W^s(p) \cap U \neq \emptyset$ for
some hyperbolic periodic point $p$. Suppose that $U$ is not contained
in the closure of $W^s(p)$, we will derive a contradiction. Let $V$ be
a connected component of $M \backslash \bar{W^s(p)}$ whose
intersection with $U$ is non-empty. Such $V$ exists by our
assumption. Let $\hat{V}$ be the prime end extension of $V$, then $U$,
as a subset of $\hat{V}$, is an open neighborhood of an arc in the
prime ends. Since the rotation number on the prime end circle is
irrational, this implies, by the the above lemma, $U$ contains
infinitely many priodic points. But $U$ is periodic point free, by our
assumption, this contraction shows that $U \subset \bar{W^s(p)}$.

This proves our main theorem.

\section{Other problems and conjectures of Poincar\'e}

Poincar\'e already noted the importance of the generic properties
area-preserving diffeomorphisms in his study of the three-body
problem. The following two fundamental conjectures are due to
Poincar\'e \cite{Po1892}.

\begin{conj} For generic $C^r$ area-preserving
diffeomorphisms on compact surface $M$, the set of all periodic points
are dense.
\end{conj}

\begin{conj} There exist a residual set $R \in \dfr$
such that if $f \in R$ and $p$ is a hyperbolic periodic point of $f$ then
the homoclinic points of $p$ is dense in both stable and unstable manifolds
of $p$. In other words, let $J$ be a segment in $W^s(p)$ (or $W^u(p)$),
then $\overline{W^s(p) \cap W^u(p) \cap J} =J$.
\end{conj}

In $C^1$ topology, both of the above conjectures are proved to be
true. The first one is a consequence of the so-called closing lemma. It
is proved by Pugh \cite{Pu67} and later improved to various cases by
Pugh \& Robinson \cite{PR83}.  A different proof was given by Liao
\cite{Liao79} and Mai \cite{Mai86}. The second conjecture in $C^1$
topology is a result of Takens \cite{Takens72}. The high dimensional
analog was proved by Xia \cite{Xi96a}. It can also be regarded as a
so-called $C^1$ connection lemma, first proved by Hayashi \cite{Ha97}
and later simplified and generalized by Xia \cite{Xi96a}, Wen \& Xia
\cite{WX99} \cite{WX00}.

In $C^r$ topology with $r > 1$, little progress has been made for
these two conjectures and it's known to be an extremely difficult
problem (cf Smale).  The local perturbation methods used in the $C^1$
case no longer seem to work and examples suggests that a more global
approach has to be developed (Gutierrez \cite{Gutierrez87}).

While unable to prove these two conjectures, our main result offers a
very strong evidence supporting the conjectures. Moreover, as an easy
consequence, our result implies that the conjecture two implies
conjecture one.


\begin{thebibliography}{10}

\bibitem{Arnold78}
V.~I. Arnold.
\newblock {\em Mathematical Methods of Classical Mechanics}.
\newblock Springer-Verlag, 1978.

\bibitem{Birkhoff66}
G.~D. Birkhoff.
\newblock {\em Dynamical Systems}, volume~9.
\newblock American Math.\ Soc.\ Colloquium Publications, 1966.

\bibitem{LY97}
P.~Le Calvez and J.C. Yoccoz.
\newblock Un th\'eor\`me d'indice pour les hom\'omorphismes du plan au
  voisinage d'un point fixe.
\newblock {\em Ann.\ of Math. (2)}, 146(2):241--293, 1997.

\bibitem{Ca13}
C.~Caratheodory.
\newblock uber die begrenzung einfach zusammenhangender gebiete.
\newblock {\em Math.\ Ann.}, 73:323--370, 1913.

\bibitem{CZ83}
C.~Conley and E.~Zehnder.
\newblock The birkhoff-lewis fixed point theorem and a conjecture of v. i.
  arnold.
\newblock {\em Invent.\ Math.}, 73(1):33--49, 1983.

\bibitem{Franks96}
J.~Franks.
\newblock Rotation vectors and fixed points of area preserving surface
  diffeomorphisms.
\newblock {\em Trans.\ Amer.\ Math.\ Soc.}, 348(7):2637--2662, 1996.

\bibitem{Franks99}
J.~Franks.
\newblock The conley index and non-existence of minimal homeomorphisms.
\newblock {\em Illinois Journal of Math.}, 43:457--464, 1999.

\bibitem{FL03}
J.~Franks and P.~Le Calvez.
\newblock Regions of instability for non-twist maps.
\newblock {\em Ergodic Theory Dynam. Systems}, 23(1):111--141, 2003.

\bibitem{Gutierrez87}
C.~Gutierrez.
\newblock A counter-example to a $c\sp 2$ closing lemma.
\newblock {\em Ergodic Theory \& Dynamical Systems}, 7(4):509--530, 1987.

\bibitem{Ha97}
S.~Hayashi.
\newblock Connecting invariant manifolds and the solution of the $c^1$
  stability and $\omega$-stability conjectures for flows.
\newblock {\em Ann. of Math.}, 145(1):81--137, 1997.

\bibitem{Liao79}
S.T.\ Liao.
\newblock An extension of the $c^1$ closing lemma.
\newblock {\em Acta Sci. Natur. Univ. Pekinensis}, 2:1--41, 1979.

\bibitem{Mai86}
J.~Mai.
\newblock A simpler proof of $c^1$ closing lemma.
\newblock {\em Scientia Sinica}, 10:1021--1031, 1986.

\bibitem{Mather82a}
J.~Mather.
\newblock Topological proofs of some purely topological consequences of
  carath\'eodory's theory of prime ends.
\newblock {\em in Selected Studies. Eds.\ Th.\ M. Rassias and G. M. Rassias},
  pages 225--255, 1982.

\bibitem{Oliveira87}
F.~Oliveira.
\newblock On the generic existence of homoclinic points.
\newblock {\em Ergod.\ Th.\ \& Dynam.\ Sys.}, 7:567--595, 1987.

\bibitem{Oliveira00}
F.~Oliveira.
\newblock On $c^\infty$ genericity of homoclinic orbits.
\newblock {\em Nonlinearity}, 13:653--662, 2000.

\bibitem{Pixton82}
D.~Pixton.
\newblock Planar homoclinic points.
\newblock {\em J. Diff.\ Equations}, 44:1365--382, 1982.

\bibitem{Po1892}
H.~Poincar\'e.
\newblock {\em Les m\'ethodes nouvelles de la m\'ecanique c\'eleste}.
\newblock Paris, 1892.

\bibitem{Pu67}
C.~Pugh.
\newblock The closing lemma.
\newblock {\em Amer. J. Math.}, 89:956--1021, 1967.

\bibitem{PR83}
C.~Pugh and C.~Robinson.
\newblock The $c^1$ closing lemma, including hamiltonians.
\newblock {\em Ergod. Th. \& Dynam. Sys.}, 3:261--313, 1983.

\bibitem{Robinson70}
C.~Robinson.
\newblock Generic properties of conservative systems, i, ii.
\newblock {\em Amer.\ J. of Math}, 92:562--603, 897--906, 1970.

\bibitem{Robinson73}
C.~Robinson.
\newblock Closing stable and unstable manifolds on the two-sphere.
\newblock {\em Proc.\ Amer.\ Math.\ Soc.}, 41:299--303, 1973.

\bibitem{SM71}
C.L. Siegel and J.K. Moser.
\newblock {\em Lectures on Celestial Mechanics}.
\newblock Springer, 1971.

\bibitem{Takens72}
F.~Takens.
\newblock Homoclinic points in conservative systems.
\newblock {\em Invent.\ Math.}, 18:267--292, 1972.

\bibitem{WX99}
L.~Wen and Z.~Xia.
\newblock A basic $c\sp 1$ perturbation theorem.
\newblock {\em J. Differential Equations}, 154(2):267--283, 1999.

\bibitem{WX00}
L.~Wen and Z.~Xia.
\newblock On $c^1$ connecting lemmas.
\newblock {\em Trans. Amer. Math. Soc.}, 352(10), 2000.

\bibitem{Xi96a}
Z.~Xia.
\newblock Homoclinic points in symplectic and volume-preserving diffeomorphism.
\newblock {\em Commun. Math. Phys.}, 177:435--449, 1996.

\end{thebibliography}
\end{document}